\documentclass[12pt,twoside]{article}
\usepackage[latin1]{inputenc}
\usepackage{times}
\usepackage{epsfig}
\usepackage{dsfont}
\usepackage{amsfonts}
\usepackage{float}
\usepackage{amsmath}
\usepackage{latexsym}
\usepackage{xr} 
\usepackage{natbib}
\usepackage{color}
\usepackage{multirow}
\usepackage{pgfplots}
\usepackage{placeins}
\usepackage{subfigure}
\usepackage{tabularx}
\usepackage{graphicx}
\usepackage{booktabs} 
\usepackage{siunitx} 
\usepackage{pgfplotstable} 
\usepackage{natbib}
\usepackage{enumerate}

\usepgfplotslibrary{polar}
\usepgflibrary{shapes.geometric}
\usetikzlibrary{calc}
\pgfplotsset{my style/.append style={axis x line=middle, axis y line=
middle, xlabel={$x$}, ylabel={$y$}, axis equal }}

\setcitestyle{square,notesep={; }}

\relax 
\allowdisplaybreaks[3]

\newtheorem{definition}{{\sc Definition}\sc}[section]
\newcommand{\bdefi}{\begin{definition}}
\newcommand{\edefi}{\end{definition}}

\newtheorem{appropr}[definition]{{\sc Approximation Procedure}\sc}
\newcommand{\bappr}{\begin{appropr}}
\newcommand{\eappr}{\end{appropr}}

\newtheorem{bedi}[definition]{{\sc Condition}\sc}
\newcommand{\bbd}{\begin{bedi}}
\newcommand{\ebd}{\end{bedi}}

\newtheorem{bedin}[definition]{{\sc Conditions}\sc}
\newcommand{\bbdn}{\begin{bedin}}
\newcommand{\ebdn}{\end{bedin}}

\newtheorem{corollary}[definition]{{\sc Corollary}\sc}
\newcommand{\bco}{\begin{corollary}}
\newcommand{\eco}{\end{corollary}}

\newtheorem{lemma}[definition]{{\sc Lemma}\sc}
\newcommand{\blem}{\begin{lemma}}
\newcommand{\elem}{\end{lemma}}

\newtheorem{proposition}[definition]{{\sc Proposition}\sc}
\newcommand{\bpro}{\begin{proposition}}
\newcommand{\epro}{\end{proposition}}

\newtheorem{satz}[definition]{{\sc Theorem}\sc}
\newcommand{\bsa}{\begin{satz}}
\newcommand{\esa}{\end{satz}}

\newtheorem{theorem}[definition]{{\sc Theorem}\sc}
\newcommand{\bth}{\begin{theorem}}
\newcommand{\eth}{\end{theorem}}

\newtheorem{assumption}[definition]{{\sc Assumption}\sc}
\newcommand{\bas}{\begin{assumption}}
\newcommand{\eas}{\end{assumption}}

\newtheorem{assumptions}[definition]{{\sc Assumptions}\sc}
\newcommand{\bass}{\begin{assumptions}}
\newcommand{\eass}{\end{assumptions}}

\newtheorem{abb}{{\sc Figure}\sc}
\newcommand{\babb}{\begin{abb}}
\newcommand{\eabb}{\end{abb}}

\newenvironment{remark}{\begin{rmk}\sl}{\end{rmk}}
\newtheorem{rmk}{{\sc Remark}\sc}[section]
\newcommand{\brem}{\begin{remark}}
\newcommand{\erem}{\end{remark}}

\newenvironment{remarks}{\begin{rmks}\sl}{\end{rmks}}
\newtheorem{rmks}{{\sc Remarks}\sc}[section]
\newcommand{\brems}{\begin{remarks}}
\newcommand{\erems}{\end{remarks}}

\newenvironment{example}{\begin{exmp}\rm}{\end{exmp}}
\newtheorem{exmp}{{\sc Example}\sc}[section]
\newcommand{\bbsp}{\begin{example}}
\newcommand{\ebsp}{\end{example}}
\newcommand{\bexa}{\begin{example}}
\newcommand{\eexa}{\end{example}}

\newtheorem{model}{{\sc Model}\sc}[section]
\newcommand{\bmdl}{\begin{model}}
\newcommand{\emdl}{\end{model}}

\newtheorem{scheme}{{\sc Scheme}\sc}[section]
\newcommand{\bscm}{\begin{scheme}}
\newcommand{\escm}{\end{scheme}}

\newenvironment{tabelle}{\begin{tabl}\rm}{\end{tabl}}
\newtheorem{tabl}{{\bf Table}}
\newcommand{\btab}{\begin{tabelle}}
\newcommand{\etab}{\end{tabelle}}

\newenvironment{exercise}{\begin{exc}\sl}{\end{exc}}
\newtheorem{exc}{Exercise}[section]
\newcommand{\bexe}{\begin{exercise}}
\newcommand{\eexe}{\end{exercise}}

\relax

\newcommand{\renu}{\mathbb{R}}

\newcommand{\qed}{\mbox{ } \hfill $\Box$\\ }

\newcommand{\bay}{\begin{array}}
\newcommand{\eay}{\end{array}}

\newcommand{\bqa}{\begin{eqnarray*}}
\newcommand{\eqa}{\end{eqnarray*}}

\newcommand{\bee}{\begin{eqnarray*}}
\newcommand{\eee}{\end{eqnarray*}}

\newcommand{\bea}{\begin{eqnarray*}}
\newcommand{\eea}{\end{eqnarray*}}

\newcommand{\bqan}{\begin{eqnarray}}
\newcommand{\eqan}{\end{eqnarray}}

\newcommand{\be}{\begin{eqnarray}}
\newcommand{\ee}{\end{eqnarray}}

\newcommand{\bit}{\begin{itemize}}
\newcommand{\eit}{\end{itemize}}

\newcommand{\ben}{\begin{enumerate}}
\newcommand{\een}{\end{enumerate}}

\newcommand{\beq}{\begin{equation}}
\newcommand{\eeq}{\end{equation}}

\newcommand{\bdes}{\begin{description}}
\newcommand{\edes}{\end{description}}

\newcommand{\btb}{\begin{tabular}}
\newcommand{\etb}{\end{tabular}}

\newcommand{\bcen}{\begin{center}}
\newcommand{\ecen}{\end{center}}

\newcommand{\bmp}{\begin{minipage}}
\newcommand{\emp}{\end{minipage}}

\newcommand{\olX}{\overline{X}}

\newcommand{\vX}{\boldsymbol{X}}

\newcommand{\vmu}{\boldsymbol{\mu}}

\newcommand{\vSigma}{\boldsymbol{\Sigma}}

\newcommand{\veins}{{\bf 1}}

\newcommand{\whsigma}{\widehat{\sigma}}

\numberwithin{equation}{section}

\begin{document}


	\title{\large \bf Permuting Incomplete Paired Data:\\A Novel Exact and Asymptotic Correct Randomization Test 
	}
	\author{Lubna Amro$^{*}$ and  Markus Pauly$^{*}$ \\[1ex] 
	}
	\maketitle

	\begin{abstract}
	Various statistical tests have been developed for testing the equality of means in matched pairs with missing values. However, most existing methods are commonly based on certain distributional assumptions such as normality, 0-symmetry or 
	homoscedasticity of the data. The aim of this paper is to develop a statistical test that is 
	robust against deviations from such assumptions and also leads to valid inference in case of heteroscedasticity or skewed distributions. 
	This is achieved by applying a novel randomization approach. The resulting test procedure is not only shown to be asymptotically correct but is also finitely exact 
	if the distribution of the data is invariant with respect to the considered randomization group. 
	Its small sample performance is further studied in an extensive simulation study and compared to existing methods. 
	Finally, an illustrative data example is analyzed. 
		
	\end{abstract}
	
	
	\vfill
	\vfill
	
	\noindent${}^{*}$ {University of Ulm, Institute of Statistics, Germany\\
		\mbox{ }\hspace{1 ex}email: markus.pauly@uni-ulm.de, lubna.amro@uni-ulm.de}

\section{Motivation and introduction}\label{int}

Matched pairs designs are concerned with scientific experiments in which subjects are observed repeatedly under two different treatments or time points. 
Testing equality of means of such data, where some of the components are randomly missing is a problem usually encountered in practice. In previous years, several statistical methods have been proposed for dealing with this issue. Here, the most simplest suggestion is to only work with the completely observed pairs and to carry out the 
well-known paired t-test. The main problem of this approach is that it discards all information from the incomplete observations. 
This will typically lead to a power loss. Moreover, the test can only be applied when the data are assumed to be missing completely at random; 
see the monograph of \cite{little2014statistical} for an introduction and discussion on the different missingness mechanisms.  
One possibility to overcome this problem might be to apply multiple-imputation techniques, see e.g. \cite{verbeke2009linear}. However these techniques often 
require large sample sizes for being correct, see e.g. the discussion in \cite{akritas2002nonparametric}. 
Other possibilities that use all observed information in the matched pairs design (but not more) include the statistical tests proposed in 
\citet{lin1974difference}, \citet{ekbohm1976oncomparing}, \citet{bhoj1978testing}, \citet{looney2003method}, \citet{kim2005statistical}, and \citet{samawi2014notes}. 
All of these tests, however, rely on specific model assumptions such as symmetry or even bivariate normality.
Especially the latter is typically not met in practice but crucial for relatively small sample sizes. 
If these specific assumptions are not met the corresponding tests might possess possibly inflated type-I error control. 

Therefore, the aim of the present paper is to provide a statistical test for 
partially paired data which does not require any parametric assumptions, possesses nice small and large sample properties and (only) uses all observed information. 
This is achieved by applying a 
novel studentized randomization approach to a weighted test statistic involving both, the paired $t$-test as well as the Welch test statistic. 
The resulting randomization test combines recent findings on studentized permutation tests for complete observations in unpaired and paired designs, see 
\cite{janssen1997}, \cite{janssen1999nonparametric}, \cite{janssen1999testing}, \cite{neubert2007}, \citet{pauly2011discussion}, \citet{konietschke2012}, 
\citet{omelka2012testing}, \cite{chung2013}, \citet{diciccio2015robust}, \cite{PBK} and \citet{chung2016asymptotically}. 
The application of the studentized randomization approach to incomplete data is original to the present paper. \\

To formulate the concrete testing problem, let us consider a general matched pairs design given by i.i.d random vectors 
\bqan \label{model}
\mathbf{X}_j= \binom{X_{1j}}{X_{2j}},\quad j=1,\ldots,n,
\eqan 
with mean vector $E(\mathbf{X}_1)=\vmu=(\mu_1,\mu_2)'\in \renu^2$ and an arbitrary covariance matrix
$Cov(\mathbf{X}_1)=\vSigma>0$. Within Model~\eqref{model} we now like to test the null hypothesis $H_0:{\{\mu_1=\mu_2\}}$. Here, some of the random vectors components might be missing completely at random as described in the next section.\\

This paper is organized as follows: In Section \ref{mod} we introduce the statistical model for incomplete paired data and a weighted test statistic. 
A permutation test based on this statistic is introduced in Section \ref{perm}, where also its theoretical properties are explained. 
We review the existing methods for incomplete paired data in Section \ref{alt} to choose adequate competitors to investigate together with our novel procedure in
extensive simulations given in Section \ref{sim}. Finally, a real data example is considered in Section \ref{example}. All proofs are shifted to the Appendix. \\

\section{Statistical model, hypotheses and statistics} \label{mod}

When some of the components of the matched pairs are missing completely at random, 
we may sort the data into complete and incomplete observed random vectors and rewrite Model~\eqref{model} with independent random variables as
\bqan \label{model: missing}
\underbrace{\binom{X_{11}^{(c)}}{X_{21}^{(c)}},\dots, \binom{X^{(c)}_{1n_1}}{X^{(c)}_{2n_1}}}_{\vX^{(c)}}, \underbrace{\binom{X_{11}^{(i)}}{--},\dots,\binom{X_{1n_2}^{(i)}}{--},\binom{--}{X_{21}^{(i)}},\dots,\binom{--}{X_{2n_3}^{(i)}}}_{\vX^{(i)}},
\eqan 
where $n=n_1+n_2+n_3$. As in Model \eqref{model} we assume that the first components $X_{1j}^{(c)},X_{1k}^{(i)}$ are i.i.d. with mean 
$\mu_1$ and variance $\sigma_1^2\in(0,\infty)$ and the second components
$X_{2j}^{(c)},X_{2\ell}^{(i)}$ are i.i.d. with mean $\mu_2$ and variance $\sigma_2^2\in(0,\infty)$ for $j=1,\dots,n_1, k=1,\dots,n_2, l=1,\dots,n_3$. 
Moreover, the complete pairs $(X_{1j}^{(c)},X_{2j}^{(c)})'$ are i.i.d. with mean vector $\vmu=(\mu_1,\mu_2)'$ and some unstructured covariance matrix $\vSigma>0$. 
Thus, the observations are missing completely at random. 
At the end of Section~\ref{perm} we comment on how to relax these model assumptions to the case where the incomplete observations are only assumed to be independent.

In the setting (\ref{model: missing}), we like to use all the available data to test the null hypotheses $H_0:{\{\mu_1=\mu_2\}}$ against the one-sided alternative ${\{\mu_1>\mu_2\}}$ 
or the two-sided alternative ${\{\mu_1\neq\mu_2\}}$. For ease of convenience, we below focus on the one-sided case and note that two-sided tests can be obtained 
similarly by, for example, taking absolute values.\\
To motivate our test statistic, we first consider the most extreme situations. In a completely observed case with $n_1=n$, the matched pairs would typically be inferred with the usual paired $t$-test type statistic
\bqan \label{T1}
T_{1}= T_1(\vX^{(c)}) = \frac{n_1^{-1}\sum_{j=1}^{n_1} D_j}{\sqrt{\whsigma^2/n_1}} =  \frac{\overline{D}_\cdot}{\sqrt{\whsigma^2/n_1}},
\eqan
where $D_j=X_{1j}^{(c)}-X_{2j}^{(c)}$ denote the differences of the first and second component for $j=1,\ldots,n_1$ and 
$\whsigma^2=(n_1-1)^{-1}\sum_{j=1}^{n_1} (D_j-\overline{D}_\cdot)^2$ is the empirical variance of these differences.\\
If, in comparison, only incomplete vectors were observed, i.e. $n=n_2+n_3$, we would be in the situation of an extended Behrens-Fisher problem and could apply a Welch-type test in the test statistic
\bqan \label{T2}
T_{2}= T_2(\vX^{(i)}) = \frac{\overline{X}_{1\cdot}^{(i)} - \overline{X}_{2\cdot}^{(i)}}{\sqrt{\whsigma_1^2/n_2 + \whsigma_2^2/n_3}},
\eqan
where $\overline{X}_{1\cdot}^{(i)} = n_2^{-1}\sum_{k=1}^{n_2} X_{1k}^{(i)}$ and $\overline{X}_{2\cdot}^{(i)} = n_3^{-1}\sum_{\ell=1}^{n_3} X_{2\ell}^{(i)}$ are the sample means and 
$\whsigma_1^2=  (n_2-1)^{-1}\sum_{k=1}^{n_2} (X_{1k}^{(i)} - \overline{X}_{1\cdot}^{(i)})^2$ and 
$\whsigma_2^2 =  (n_3-1)^{-1}\sum_{\ell=1}^{n_3} (X_{2\ell}^{(i)} - \overline{X}_{2\cdot}^{(i)})^2$ are the corresponding empirical variances. To deal with all other cases of \eqref{model: missing} adequately, it is then
reasonable to propose the usage of a weighted test statistic
\bqan \label{T}
T = T(\vX^{(c)},\vX^{(i)}) = \sqrt{a} T_{1}(\vX^{(c)}) + \sqrt{1-a} T_2(\vX^{(i)}),
\eqan
where $a\in[0,1]$ can be used to weight the inference drawn from the complete pairs. 
We suggest to use $a = 2n_1/(n+n_1)$ or $a=n_1/n$ to cover the usually applied methods in the two extreme cases $n_1=n$ and $n_1=0$. 
The asymptotic null distribution of $T$ as $\min(n_1,n_2,n_3)\to \infty$ is stated below. The asymptotic framework is, e.g., reasonable when assuming that each component is missing at random with a certain probability $p\in(0,1)$.

\begin{theorem}\label{theo: clt T}
	Suppose that $n_2/(n_2+n_3)\to \kappa\in(0,1)$ as $\min(n_1,n_2,n_3)\to \infty$ then we have under $H_0: \mu_1=\mu_2$
	\bqan
	\sup_{x \in \renu} \left| P(T\leq x) - \Phi(x) \right| \to 0,
	\eqan
	where $\Phi$ denotes the cdf of $N(0,1)$.
\end{theorem}

From this a consistent asymptotic level $\alpha$ test is given by $\varphi = \veins\{T> z_{1-\alpha}\}$ for testing the one-sided case  ${\{\mu_1\leq\mu_2\}}$ vs ${\{\mu_1>\mu_2\}}$ and similar in the two-sided case, see also \citet{samawi2014notes}. 
Here $z_{1-\alpha}$ denotes the $(1-\alpha)$-quantile of the standard normal distribution, $\alpha\in(0,1)$.
However, large sample sizes of the incomplete paired data are necessary to get an accurate approximation of this test statistic. In particular, if the sample size is small or  the underlying 
distribution shows deviation  from symmetry, the test will in general not keep the preassigned type-I error level accurately, see the simulation study in Section \ref{sim} below. Therefore, we suggest a novel permutation
approach in the next section that is not only asymptotically valid but also finitely exact under certain assumptions.

\section{The permutation procedures}\label{perm}

Previous to our work, permutation techniques have been applied to the incomplete pairs design by \citet{maritz1995permutation} and \citet{yu2012}. 
However, these methods have the drawback that certain distributional assumptions such as 0-symmetry and equal variances or sample sizes are needed to obtain 
a (at least asymptotically) valid level $\alpha$ test.\\

To overcome this issue, we propose a different permutation test that has the finite sample exactness property under similar invariance assumptions but is also (asymptotically) robust against deviations such as heteroscedasticity or underlying skewed distributions, 
i.e. valid for our general model~\eqref{model: missing}. The procedure is based on separate results about studentized permutation tests for the paired and  unpaired case. For example, 
in the classical matched pairs case with $n_1=n$ \citet{janssen1999testing} and \citet{KP14} have already recommended (based on theoretical results and extensive simulation studies) a studentized permutation test in the paired $t$-test, 
where the components of each pair are randomly permuted. In comparison, in the general Behrens-Fisher set-up with $n_2+n_3=n$, \citet{janssen1997, janssen2003} as well as \citet{janssen2005} proposed a studentized permutation test in 
the Welch-type statistic $T_2$ that is based on randomly permuting the pooled sample, see also \citet{chung2013} as well as \citet{PBK} for recent generalizations of this approach to multiple samples and general factorial designs. 
These studentized permutation procedures have also been implemented in R and SAS, see \citet{placzekstudentisierte} as well as the R package GFD. \\

Jointly applying the ideas of both approaches we propose a permutation test that is based on the weighted test statistic $T$ as follows: 
First, we randomly permute the components of each complete observation resulting in 
\bqan
\vX^{(c)}_\tau = \binom{X_{\tau_1(1)1}^{(c)}}{X_{\tau_1(2)1}^{(c)}},\dots, \binom{X^{(c)}_{\tau_{n_1}(1)n_1}}{X^{(c)}_{\tau_{n_1}(2)n_1}},
\eqan
where $\tau_i, i=1,\dots,n_1,$ are uniformly distributed on the symmetric group $\mathcal{S}_2$ (i.e. are random permutations of $(1,2)$) and independent of $\vX=(\vX^{(c)},\vX^{(i)})$. 
Second, for the incomplete observations, we randomly permute the pooled incomplete sample $(Z_{1},\dots,Z_{n_2+n_3})= (X_{11}^{(i)},\dots,X_{1n_2}^{(i)},X_{21}^{(i)},\dots,X_{2n_3}^{(i)})$ 
by means of a random permutation $\pi$ that is uniformly distributed on $\mathcal{S}_{n_2+n_3}$ and independent of all other variables 
resulting in 
$$
\vX^{(i)}_\pi =  (Z_{\pi(1)},\dots,Z_{\pi(n_2+n_3)}).
$$ 
Finally, we approximate the null distribution of $T$ by the conditional distribution of the permutation version of the test statistic 
$T_p= T(\vX^{(c)}_\tau,\vX^{(i)}_\pi) = \sqrt{a} T_{1}(\vX^{(c)}_\tau) + \sqrt{1-a} T_2(\vX^{(i)}_\pi)$. The validity of this approach is proven below, where '$\stackrel{p}{\longrightarrow}$' denotes convergence in probability. 

\begin{theorem}\label{theo: cclt T}
	Suppose that $n_2/(n_2+n_3)\to \kappa\in(0,1)$ as $\min(n_1,n_2,n_3)\to \infty$ then we have 
	\bqan\label{cclt}
	\sup_{x \in \renu} \left| P(T_p\leq x|\vX) - \Phi(x) \right| \stackrel{p}{\longrightarrow} 0.
	\eqan
	Moreover, choosing $c_p(\alpha)$ as $(1-\alpha)$-quantile of the permutation distribution of $T_p$ given the data, the test 
	$\varphi_p = \veins\{T> c_p(\alpha)\} + \gamma_p \veins\{T = c_p(\alpha)\}$ is consistent and of asymptotic level $\alpha$ for testing $H_0$ against one-sided alternatives ${\{\mu_1>\mu_2\}}$. Moreover, in the special case that $D_1$ is $0$-symmetric and  $X_{11}^{(i)}$ has the same distribution as $X_{21}^{(i)}$ the test is even finitely exact.
\end{theorem}

By inverting the above one-sided test, we can also construct one--sided confidence interval for the mean difference with asymptotic coverage of probability $(1-\alpha)$. Moreover, a related 
two-sided confidence interval for $\mu_1-\mu_2$ based on this approach is given by
$$
\left[ (T \pm c_P(\alpha/2))/S_w \right],
$$
where 
$$
S_w = \frac{\sqrt{a}}{\widehat{\sigma}/\sqrt{n_1}} + \frac{\sqrt{1-a}}{(\widehat{\sigma}_1/\sqrt{n_2} + \widehat{\sigma}_2/\sqrt{n_3} )}
$$
is a weighted sum of the inverses of the variance estimates for the paired and unpaired samples, respectively.
Finally, applying the intersection-union principle (see e.g. \cite{perlman1999emperor}) it is straightforward to obtain a permutation TOST test for the bioequivalence hypothesis $H_0^{(b)}: |\mu_1-\mu_2|\geq \epsilon$ 
for some given $\epsilon>0$ from the above result.

\begin{remark}[Relaxing the Model Assumptions]\label{rem1}
	If the data is only missing at random, the i.i.d. assumptions of the incompletely observed first and second components, respectively, in model (\ref{model: missing}) 
	might be violated. However, the above permutation approach can still be valid in such a situation. 
	To accept this consider the relaxed model, where the incomplete components $(X_{1k}^{(i)})_{k}$ and $(X_{2\ell}^{(i)})\ell$ 
	are independent of $\vX^{(c)}$ and assumed to be (row-wise) 
	independent and infinitesimal with mean $\mu_1$ (first component) and $\mu_2$ (second component), respectively. If both means 
	$\sqrt{n_2} (\overline{X}_{1\cdot}^{(i)}-\mu_1)$ and $\sqrt{n_3} (\overline{X}_{2\cdot}^{(i)}-\mu_2)$ are asymptotically normal with positive asymptotic variance 
	it then follows from results of \cite{janssen2003} and \cite{janssen2005} that the convergences (\ref{cclt}) and $T\stackrel{d}{\rightarrow} N(0,1)$ still hold under the null hypothesis. Thus, the proposed randomization procedure remains an asymptotic exact level $\alpha$ test. 
	Moreover, skipping the assumption of mean equality, i.e. $E({X}_{sr}^{(i)})\neq E({X}_{s1}^{(c)})$ may hold for $s=1,2$, $r=1,\dots,n_{s+1}$, 
	the procedure is even valid for testing the null hypothesis $H_0:\{E({X}_{11}^{(c)})=E({X}_{21}^{(c)}) \text{ and } E(\overline{X}_{1\cdot}^{(i)})=
	E(\overline{X}_{2\cdot}^{(i)})\}$.
	
\end{remark}

\section{Alternative approaches} \label{alt}

In this subsection, we briefly review the existing literature to select possible competitors for our subsequent simulation study. 
Several methods that (only) deal with all available data in the incomplete matched pairs design have been suggested. 
However, most of the available procedures  are based on specific  assumptions such as normality, 0-symmetry or homoscedasticity of the data. We can summarize the most commonly used methods as follows: 
\begin{enumerate}[(a)]
	\item Tests based on modified maximum likelihood estimators (cf., e.g., \citet{lin1974difference}, \cite{ekbohm1976oncomparing, ekbohm1981testing}, \citet{woolson1976monte}, or \citet{hamdan1978test}).
	\item Tests based on simple mean difference estimator (see, e.g., \citet{lin1974difference}, \citet{ekbohm1976oncomparing}, \citet{bhoj1989comparing}, \citet{looney2003method}, or \citet{uddin2015testing}).
	\item Weighted linear and nonlinear combination tests (cf., e.g., \cite{bhoj1978testing, bhoj1984difference, bhoj1991testing}, \citet{kim2005statistical}, or \citet{samawi2014notes}).
	\item Resampling and more specific randomization methods  (see, e.g., \citet{maritz1995permutation} and \citet{yu2012}).
\end{enumerate}

For ease of presentation, it is not possible to compare our suggested approach with all these procedures. Also, none of them is free from distributional assumptions and at the same time robust against deviations such as heteroscedasticity and skewed distributions. 
Therefore, we selected a few procedures for the comparisons in our simulation study below. 
First, the test statistic proposed by \citet{lin1974difference} has been chosen. It is a studentized mean difference statistic given by 
\bqa
T_{LS}=\frac{\olX_{1.}^{(c,i)}-\olX_{2.}^{(c,i)}}{\sqrt{\frac{1}{n_2+n_1}+\frac{1}{n_3+n_1}-\frac{2n_1r}{(n_2+n_1)(n_3+n_1)}}\sqrt{\frac{S_1^2+S_2^2}{n-2}}}.
\eqa
Here, $\olX^{(c,i)}_{1\cdot}$ and $\olX^{(c,i)}_{2\cdot}$ denote the means of the paired and unpaired observations for sample 1 and sample 2 respectively, $r$ is the 
empirical correlation coefficient of the completely observed pairs and $S_1^2=\sum_{j=1}^{n_1}{(X_{1j}^{(c)}-\olX_{1.}^{(c,i)})^2}+\sum_{k=1}^{n_2}{(X_{1k}^{(i)}-\olX_{1.}^{(c,i)})^2}$ 
as well as $S_2^2=\sum_{l=1}^{n_3}{(X_{2l}^{(i)}-\olX_{2.}^{(i)})^2}$ are sum of squares. Lin and Stivers motivated this statistic 
assuming normality of the data, and thus proposed to approximate the distribution of $T_{LS}$ by a $t_{n-4}$-distribution with $n-4$ degrees of freedom. 
In extensive simulation studies (e.g. \citet{dunu1994comparing}, and \cite{bhoj1989comparing, bhoj1991testing}), it turned out 
that the $T_{LS}$ maintained the type-I error rate quite accurately under homoscedasticity. Moreover, it was found to possess a much better power behaviour in 
case of data with low correlation than the following competitors: The paired t-test only depending on the complete observations, 
the independent two-sample t-test only depending on the incomplete cases,  the test by \citep{hamdan1978test} based on a modified maximum likelihood estimator 
and the weighted convex combination test given in \citep{bhoj1991testing}. \\

From the group of weighted linear and nonlinear combination tests, we selected the recently developed test by \citet{kim2005statistical}
as second competitor. The test statistic is given by the following combination of paired and unpaired mean group differences:
\bqa
t_3=\frac{n_1(\olX_{1.}^{(c)}-\olX_{2.}^{(c)})+n_H(\olX_{1.}^{(i)}-\olX_{2.}^{(i)})}{\sqrt{n_1\whsigma^2+n_H^2(\whsigma_1^2/n_2+\whsigma_2^2/n_3)}},
\eqa 
where $n_H$ is the harmonic mean of $n_2$ and $n_3$. In their paper the authors proposed to approximate the null distribution of $t_3$ with a standard normal distribution and applied 
it to data from a colorectal cancer microarray study 
for identifying differentially expressed genes, see also \citet{yu2012}.\\

Finally, we selected the asymptotic test $T$ of Section 2 as it was recommended by a simulation study in \citet{samawi2014notes} in favor of their suggested 
pooled t-test and the corrected $z$-test \citep{looney2003method}. 
A similar conclusion was drawn by \citet{guo2015comparative} who also found out in an extensive simulation study that $T$ controls the type-I error rate 
quite accurate in case of moderate sample sizes.\\ 

Other tests as the corrected $z$-test by \citet{looney2003method} or the pooled t-test by \citet{samawi2014notes} have not been included since they also showed worse performances for small sample sizes in other simulation studies, 
see e.g. \citet{samawi2014notes}, \citet{guo2015comparative} and \citet{uddin2015testing}.
We thus restricted our simulation study to comparing the small sample behavior of our suggested test with the three competitors from  above.

\begin{figure}
	\begin{center}
		\subfigure {\includegraphics[width=.45\textwidth]{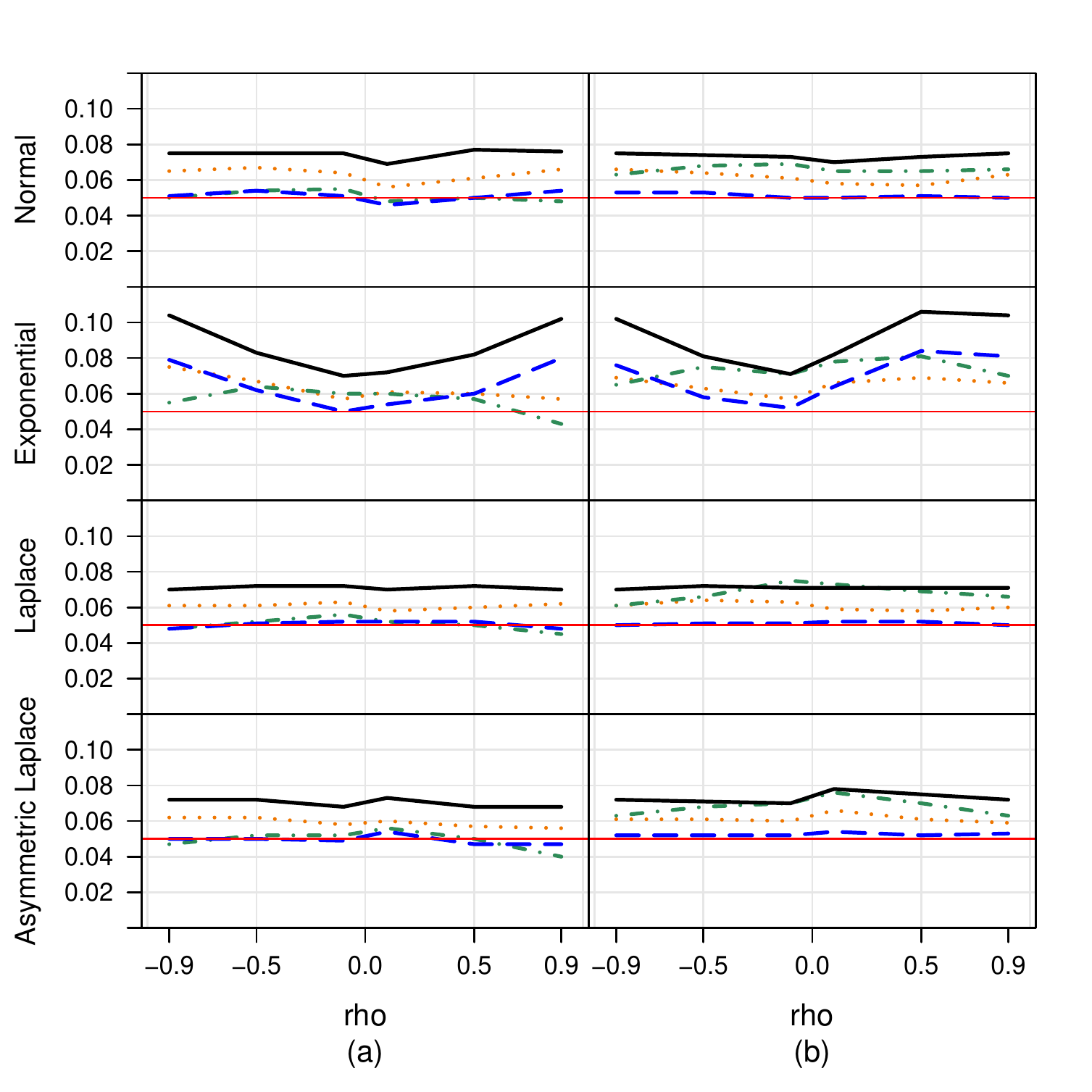}}
		\subfigure {\includegraphics[width=.45\textwidth]{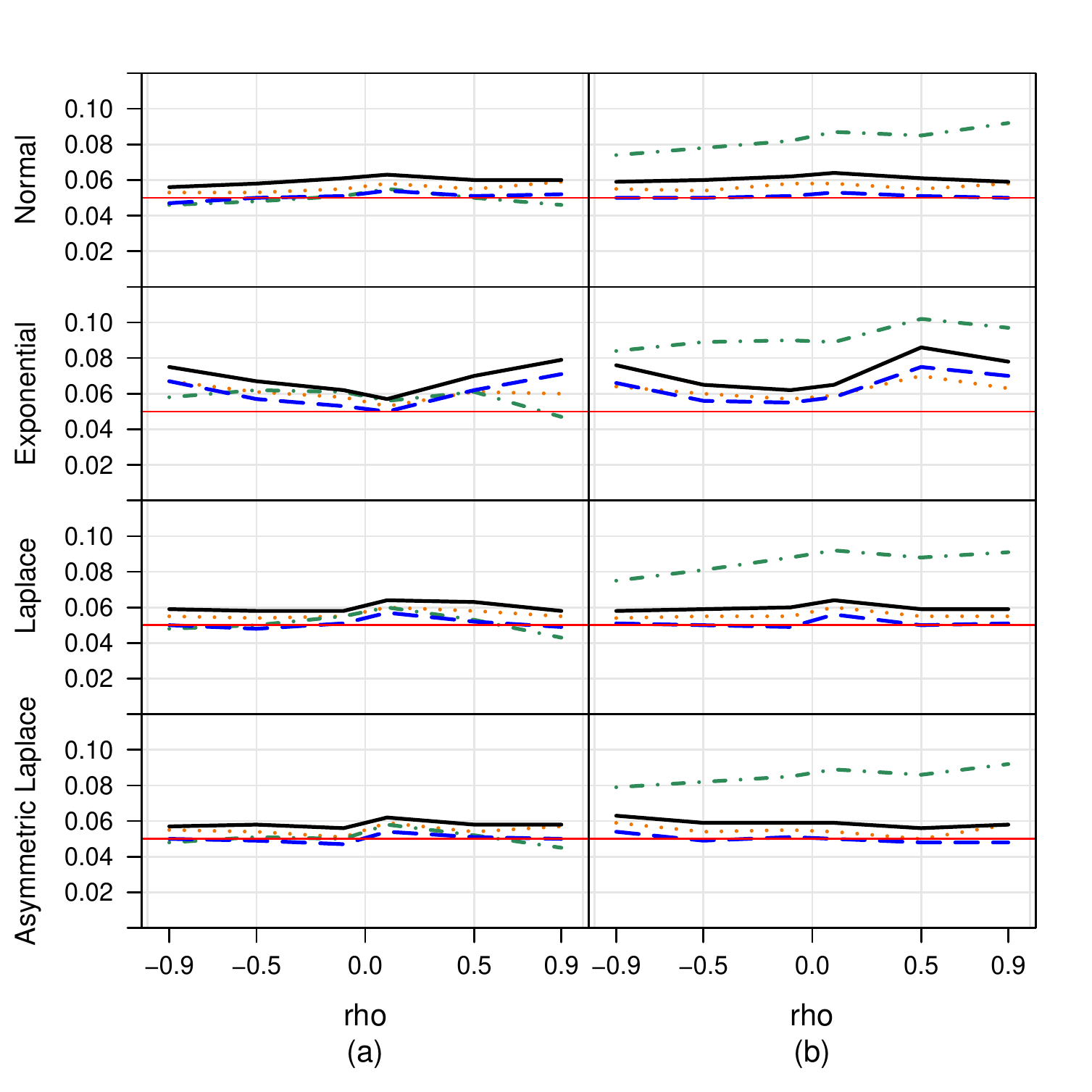}}
		\hfill
		\caption{Simulation results for the type-I error level for $T_P$ $(--)$, $T$ $(\textendash\textendash\textendash)$, 
			$T_{LS}$ $(-\cdot-)$, and $t_3$ $(\cdots)$ for different distributions under varying correlation values each with sample sizes 
			$(n_1,n_2,n_3)\in(10,10,10)$ (left) and $(30,10,10)$ (right) and different covariance matrices  $\Sigma_1$ (a) and $\Sigma_2$ (b).}
		\label{fig:multi}
	\end{center}	
\end{figure}
\section{Monte-Carlo simulations} \label{sim}
In this Section, we investigate the finite sample behavior of our suggested permutation test 
and the chosen competing procedures as described in Section \ref{alt} in extensive simulations. 
As major assessment criteria all procedures were studied with respect to (i) maintaining the preassigned type-I error level under the null hypothesis 
and (ii) their powers to detect certain alternatives. The observations were generated using four symmetric and skewed bivariate distributions: 
normal, exponential, asymmetric Laplace, and Laplace distribution each with correlation $\rho \in (-1,1)$ and with the following two covariance matrices
\begin{center}
	$\Sigma_1 = \begin{pmatrix} 1&\rho\\ \rho&1 \end{pmatrix}$,   \hspace{2pt}
	$\Sigma_2 = \begin{pmatrix} 1&\sqrt{2}\rho\\ \sqrt{2}\rho& 2\end{pmatrix}$,
\end{center}
representing a homoscedastic and a heteroscedastic setting, respectively.
The complete pairs 
were simulated from the model 
\begin{center} 
	$\begin{pmatrix}
	X^{(c)}_{1j} \\ X^{(c)}_{2j} 
	\end{pmatrix}= \Sigma^{1/2} \begin{pmatrix}\varepsilon_{1j} \\\varepsilon_{2j} \end{pmatrix} +  \begin{pmatrix}
	\mu_1\\ \mu_2
	\end{pmatrix}, \quad j=1,\dots,n_1.$
\end{center}
Here, the standardized random errors $\varepsilon_{ij}, i=1,2,$ were generated by $\varepsilon_{ij}=\frac{Y_{ij}-E(Y_{i1})}{\sqrt{var(Y_{i1})}}$, 
where $Y_{ij}$ denote normal, exponential, asymmetric Laplace, or Laplace random variables, respectively. Similarly, the incomplete observations were generated as $X^{(i)}_{2j} = \sigma_i\tilde{\varepsilon}_{ij} + \mu_i$, where $\tilde{\varepsilon}_{ij}$ 
are independent copies of the $\varepsilon_{ij}$ and the scaling factors are $\sigma_1=1$ and $\sigma_2\in\{1,\sqrt{2}\}$ depending on the choice of covariance matrix. Moreover, the following sample sizes $(n_1,n_2,n_3) \in \{(10,10,10), (30,10,10)\}$
were considered and the nominal significance level was taken as $\alpha=0.05$. Under the null hypothesis, we set $\mu_1=\mu_2=0$. To also evaluate the empirical power, the second samples were shifted by adding 
a shift $\mu_2=\delta\in\{0,0.5,1\}$. The choice of $\delta=0$ corresponds to the analysis of the empirical level of significance. 
All simulations were operated by means of the R computing environment, version 3.2.2  \citep{Rlanguage} and each setting was based on $10,000$ simulation runs and $B=1,000$ permutations runs. 
The algorithm for the computation of the $p$-value of our suggested permutation test is as follows:
\begin{enumerate}
	\item For the given partially paired data, calculate the observed test statistic $T$ with weight $a=2n_1/(n+n_1)$. 
	\item Randomly permute the components of the complete pairs $\vX^{(c)}$, resulting in $\vX^{(c)}_\tau$.
	\item Randomly permute the pooled incomplete pairs  $\vX^{(i)}$, resulting  in $\vX^{(i)}_\pi$.
	\item Calculate the value of the weighted test statistic for the permuted data $T_p=T(\vX_\tau^{(c)},\vX_\pi^{(i)})$. 
	\item Repeat the Steps $2$ and $3$ Independently $B=1,000$ times and collect the observed test statistic values in ${T_{b}, b=1,.....,B}$.
	\item Finally, estimate the two sided permutation $p$-value as follows  $p-value=min(2p_1,2-2p_1)$, where $p_1=\frac{\sum_{b=1}^{B}I(T_{pb}>= T)}{B}$. 
\end{enumerate}
\begin{figure}
	\begin{center}
		\subfigure		{\includegraphics[width=.45\textwidth]{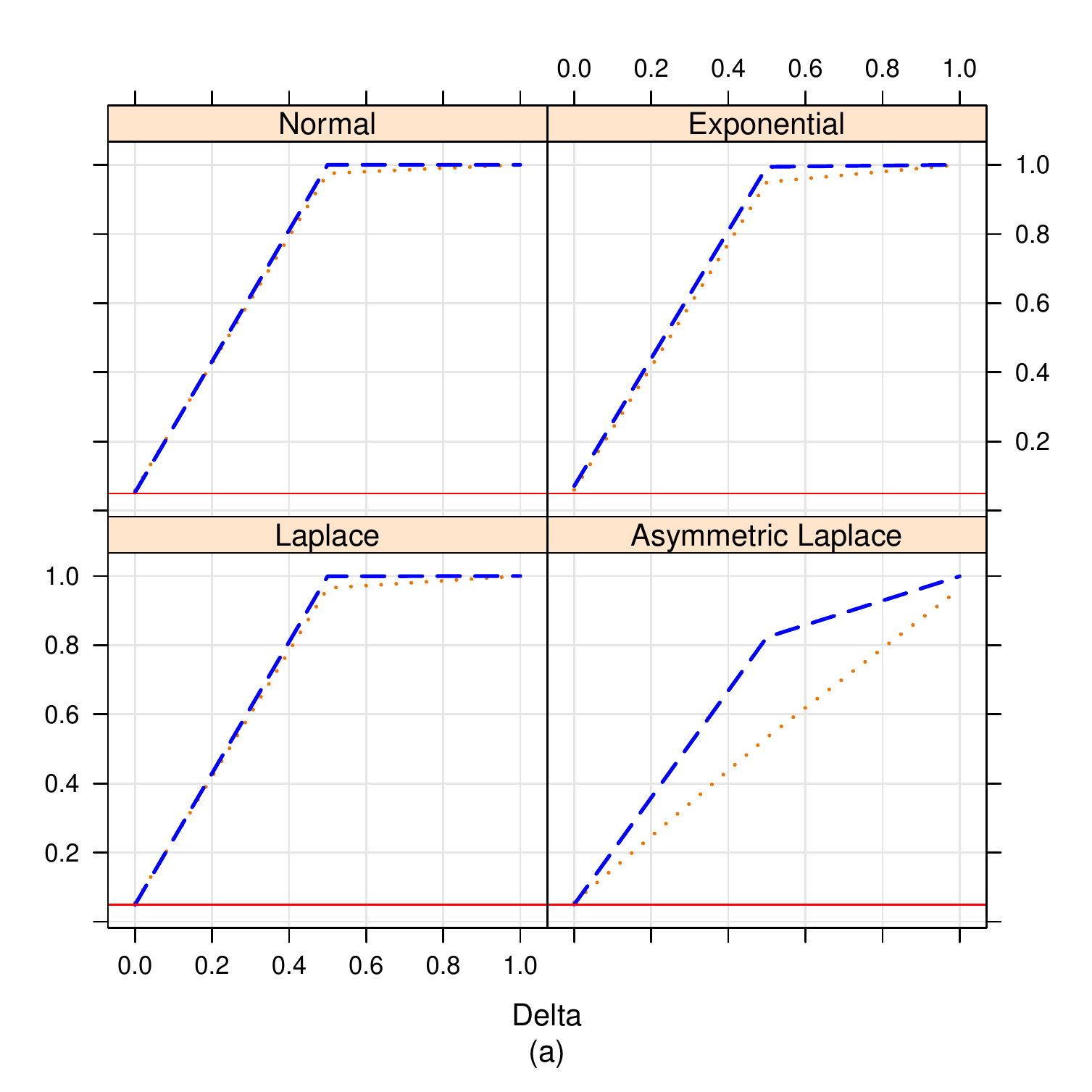}}
		\subfigure	{\includegraphics[width=.45\textwidth]{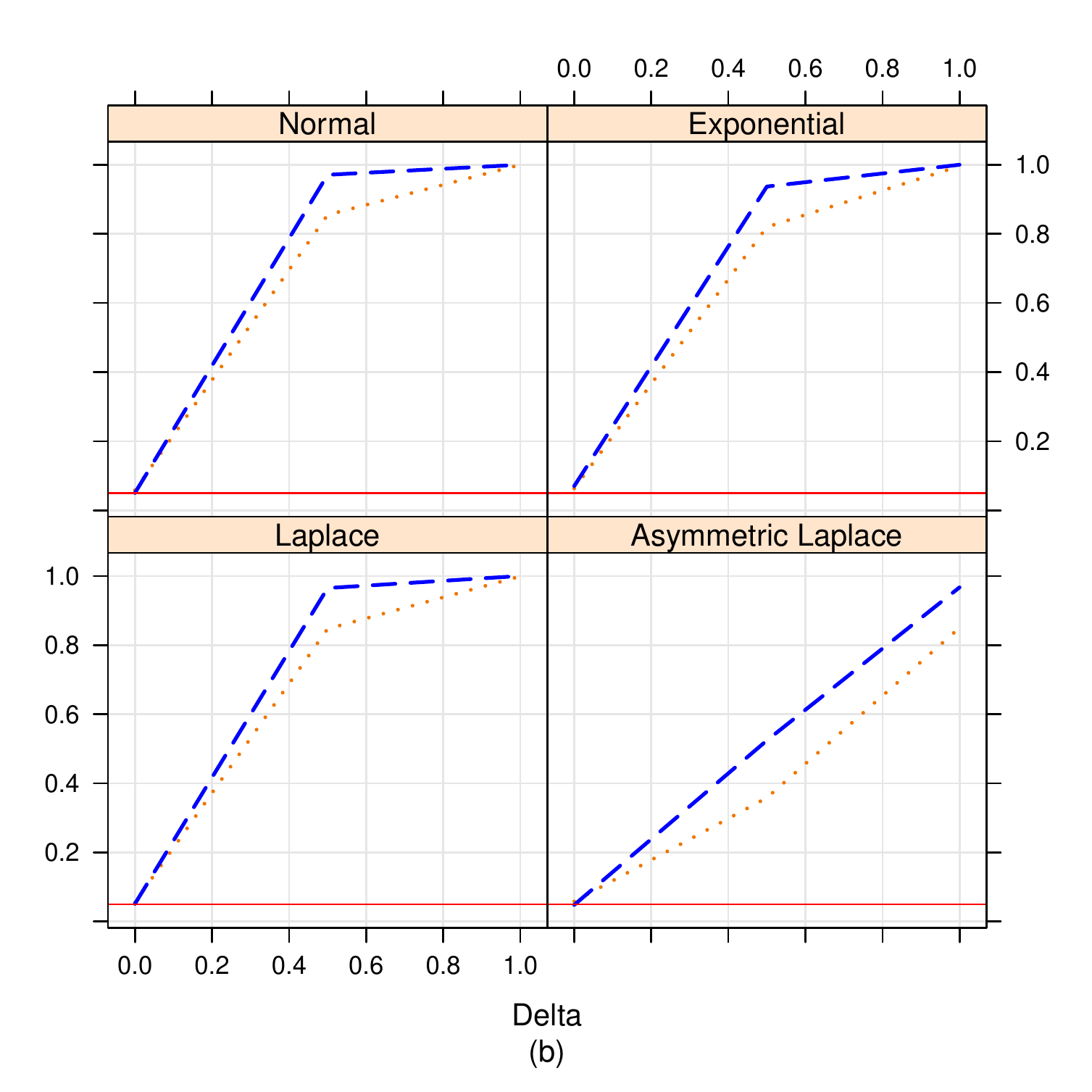}}
		\hfill
		\caption{Power simulation results $(\alpha=0.05)$ for $T_P$ $(--)$ and $t_3$ $(\cdots)$ for Symmetric and Asymmetric distributions under 
			$\rho=0.9$ for $(n_1,n_2,n_3)\in(30,10,10)$ and for  $\Sigma_1$ (a) and  $\Sigma_2$ (b).}
		\label{fig:power}
	\end{center}	
\end{figure}

The type-I error level of the proposed procedures for different sample sizes and for homoscedastic as well as heteroscedastic settings are displayed in
Figure \ref{fig:multi}. It can be readily seen from Figure \ref{fig:multi} that for all the considered settings, our suggested studentized permutation test $T_p$ tends to
result in an adequate exact level $\alpha$ test under homoscedasticity as well as heteroscedasticity and 
over the range of correlations $\rho$ for most settings. Only in the case of the skewed exponential distribution, the control is not adequate. However, 
in this case all the other competing procedures also failed to control the type-I error rate for the underlying sample sizes. 
More specific, in the case of homoscedasticity, data with low correlation ($\mathopen|\rho\mathclose|\leq0.5$), 
and exponential distribution our suggested $T_p$ still results in an accurate test decision and also performs better when
the number of complete pairs is increased.\\

In contrast, the competing tests do not control the type-I error level constantly under heteroscedasticity or even under homoscedasticity 
in the two considered cases of sample sizes. It can also be seen from Figure \ref{fig:multi} that the weighted test $T$ is rather liberal in the case of 
small number of complete pairs and less liberal for larger number of complete pairs. 
This behavior of the test does not depend on the homoscedasticity assumption. Moreover, the test that is based on the simple mean difference estimator 
$T_{LS}$ controls type-I rate in the case of equal variances under all the simulation settings except for the skewed exponential distribution. 
However, it is very liberal in the case of heteroscedasticity and its behavior is much worse for larger number of complete pairs. 
Finally, the test $t_3$ tends to result in a liberal decision in the case of small number of complete pairs and vice versa it is nearly exact level $\alpha$ 
in case of a larger number of complete pairs; especially under homoscedasticity.\\

The power of the four tests has also been studied for all settings. Due to the rather liberal behavior of the test based on $T_{LS}$ and the asymptotic test 
based on $T$, their power functions are not really comparable to the others. 
Therefore, we restrict our power simulation study to the analysis of $T_p$ and $t_3$. Moreover, since their type-I error rate control is rather similar 
for $(n_1,n_2,n_3)\in(30,10,10)$ we only consider the unbalanced case. 
Figure \ref{fig:power} shows the power of both methods for highly correlated data involving homoscedastic as well as heteroscedastic settings. 
It can be easily seen that our suggested studentized permutation test $T_p$ has a larger power under all studied settings. 
However, both tests yielded an almost similar power behavior when equal variances and symmetric distributions are present. 
In all other cases (heteroscedasticity or asymmetric distributions) the power of $t_3$ is considerably smaller. 
Note, that the power functions are not really comparable under the exponential distribution due to their rather liberal behavior and we have only included it 
for illustration and completeness. In the power simulations for all the other considered situations (not shown here) both tests performed nearly similar.

\section{Analysis of the data example}\label{example}
The Karnofsky Performance Status Scale (KPS) as proposed by \citet{karnofsky1949clinical}, is an important assessment tool for functional impairment.
This scale can  be used  for determining the prognosis of individual patients or for comparing effectiveness of several therapies.
The score ranges from 0 to 100. The lower the KPS score, the lower the probability of survival.\\
A clinical  study by \citet{hermann2001effectiveness} has been performed to assess the functional  status of hospice patients during their last seven days of life and to discover and examine the different interventions 
that have been implemented to treat the patients symptoms. Their study consists of 100 hospice patients. 
KPS scores were used to assess the functional status of the patients on the day before they died and their last day of life resulting in paired data. These pairs were not complete. For example, there were a total of $n_1=9$ complete pairs and two 
unpaired samples for the patients who were only observed on  their  day  before they died or their last  day of life  of sizes $n_2=28$ and $n_3=23$ respectively. The data was published in \citet{rempala2006asymptotic}. 
We applied all considered testing methods $T_p$, $T$, $T_{LS}$, and $t_3$ to  detect the null hypothesis of equality of means of the KPS scores on the day before they died and on the last day of their life ($H_0:{\{\mu_1=\mu_2\}}$) against 
the one-sided alternative ${\{\mu_1>\mu_2\}}$ or the two-sided alternative ${\{\mu_1\neq\mu_2\}}$. The results are summarized in Table \ref{pvalue}. 

\begin{table}[!ht]
	\caption{P-values of  the KPS Data} 
	\centering 
	\begin{tabular}{l c c } 
		\hline\hline 
		Method & One-sided p-value & Two-sided p-value \\ [0.5ex] 
		\hline \hline 
		$T_p$ & 0.0110 & 0.0220\\ 
		$T$ & 0.0025 & 0.0050\\
		$T_{LS}$ & 0.0077 & 0.0153\\
		$t_3$ & 0.0069 & 0.0138 \\
		\hline \hline 
	\end{tabular}
	\label{pvalue} 
\end{table}
All four tests provided in Table \ref{pvalue} reject the null hypothesis and indicate that the mean values of KPS score on the last day and on the day before 
they died are significantly different (two sided p-value $<0.05$). In particular, their KPS score on the day before they died was significantly greater than their score on the last day of life (right-sided p-value $<0.05$).

\section{Conclusions and discussion} \label{dis}
The problem of incomplete paired data occurs frequently in practice. In this work, we have provided a novel randomization test that is not based on any parametric 
assumptions and uses all the observed data of the matched pairs design. It was shown to be asymptotic valid and even finitely exact if certain invariance properties are met. 
For the asymptotic considerations, ideas  
of \cite{janssen1997}, \cite{janssen1999nonparametric}, \cite{janssen1999testing}, \cite{neubert2007}, \citet{pauly2011discussion}, \citet{konietschke2012}, 
\citet{omelka2012testing}, \cite{chung2013}, \citet{diciccio2015robust}, \cite{PBK} and \citet{chung2016asymptotically} 
on studentized permutation tests for paired and unpaired two-sample problems were combined. 
In addition, the theoretical asymptotic findings were also applied to bioequivalence problems and to construct one- and two-sided confidence interval 
for the two sample mean difference. 

In an extensive simulation study, the advantageous theoretical properties of the randomization test were underpinned.
There it was seen that our novel test possess nice small sample properties and particularly 
improves the small sample behavior of the asymptotic test based on the same weighted test statistic $T$.
In particular, the simulations also indicated that our randomization test is robust in terms of error control against heteroscedasticity or skewed distributions in most considered situations. 
Only in the case of the exponential distribution with rather large correlations its behaviour was rather liberal for small sample sizes. 
In this setting, however, the other considered competing procedures (\citet{samawi2014notes}, \citet{kim2005statistical}, and \citet{lin1974difference})
also failed to control the type-I error rate.

To sum up, its simulation results together with the finite exactness property under invariance and its asymptotic correctness make the new permutation procedure 
recommendable in practice.

\section*{Acknowledgement} 
This work was supported by the German Academic Exchange Service (DAAD) under the project: Research Grants - Doctoral Programmes in Germany, 2015/16 (No. 57129429).


\section*{Appendix}

\textit{Proof of Theorem~\ref{theo: clt T}:}\\ It follows from the central limit theorem and Slutzky that we have convergence in distribution
$$
T_1(\vX^{(c)}) \stackrel{d}{\longrightarrow} N(0,1) \quad\text{as well as}\quad T_2(\vX^{(i)}) \stackrel{d}{\longrightarrow} N(0,1).
$$
Since $T_1(\vX^{(c)})$ and $T_2(\vX^{(i)})$ are independent, it follows that
$$
\sqrt{a}T_1(\vX^{(c)}) + \sqrt{1-a}T_2(\vX^{(i)}) \stackrel{d}{\longrightarrow} N(0,a + 1-a) = N(0,1).
$$
The stated result thus follows from Polya's Theorem. \hfill\qed
\\

\noindent\textit{Proof of Theorem~\ref{theo: cclt T}:}\\
From the results given in \citet{janssen1999testing} and \citet{KP14} it follows under the given assumptions that 
\bqa
\sup_{x \in \renu} \left| P(T_{1}(\vX^{(c)}_\tau)\leq x|\vX^{(c)}) - \Phi(x) \right| \stackrel{p}{\to} 0.
\eqa
Similarly, the resuls from \citep{janssen2005}, see also \citep{PBK}, show that 
\bqan\label{cclt2}
\sup_{x \in \renu} \left| P(T_{2}(\vX^{(i)}_\pi)\leq x|\vX^{(i)}) - \Phi(x) \right| \stackrel{p}{\to} 0.
\eqan
Since $\pi$ and $\tau_j, j=1,\dots,n_1,$ are independent the weak convergence result follows from a conditional Slutzky Theorem, see e.g. \citet{chung2013} or \citet{pauly2009}. 
From this an application of Lemma~1 in \citet{janssen2003} proves the consistency and asymptotic exactness of the permutation test. 
The finite exactness under the given invariance assumptions follows from classical theory on randomization tests, see e.g. \citet{lehmann2006testing} or \citet{janssen2007}.
{ }\hfill\qed
\\

\noindent\textit{Proof of Remark~\ref{rem1}:}\\
First, $T\stackrel{d}{\rightarrow} N(0,1)$ follows from Raikov's Theorem. Thus, to accept the stated results we only have to argue that (\ref{cclt2}) also holds for the more general model under the stated assumptions. But this follows from Theorem~4.1 in 
\citet{janssen2005} by applying Lemma~7 in \citet{janssen2003}, see also the comments in Example~5(c) of the first paper.{ }\hfill\qed

\bibliographystyle{plainnat}
\bibliography{References}

\end{document}